\documentclass[12pt]{article}
\usepackage{amssymb}
\usepackage{amsfonts}
\usepackage{mathrsfs}
\usepackage{amsmath}
\usepackage{pifont}
\usepackage{amsopn}
\usepackage{bm}
\usepackage{amscd}
\usepackage{framed}
\usepackage{delarray}
\usepackage{xypic}
\usepackage{exscale}
\usepackage{amsthm}
\usepackage{hyperref}
\usepackage[square,comma]{natbib}

\begin{document}

\newcommand{\Proof}{\noindent {\bf Proof.}}
\newcommand{\Section}[1]{
   \stepcounter{section}
   \bigskip\noindent
   {\bf\hbox{\thesection.~}#1}\par
   \nopagebreak
   \medskip
   \renewcommand{\theequation}{\thesection.\arabic{equation}}
   \setcounter{equation}{0}
   \setcounter{subsection}{0}}
\theoremstyle{definition}
\newtheorem{definition}{Definition}[section]
\newtheorem{remark}{Remark}[section]
\theoremstyle{plain}
\newtheorem{lemma}{Lemma}[section]
\newtheorem{theorem}{Theorem}[section]
\newtheorem{corollary}{Corollary}[section]
\newtheorem{proposition}{Proposition}[section]
\newtheorem{example}{Example}[section]
\title{{\bf  A New Logic for Uncertainty}
\thanks{\it Project supported by NSFC (10731050) and PCSIRT (IRTO0742).}}
\author{\footnotesize {\begin{minipage}{5cm}\begin{center}{LUO Maokang\\Institute of Mathematics\\Sichuan University\\
Chengdu, 610064\\P.R.China}\end{center}\end{minipage}}
\begin{minipage}{1cm}{and\\ \\ \\ \\ }\end{minipage} {\begin{minipage}{5cm}\begin{center}{HE Wei\thanks{\it Corresponding author: weihe@njnu.edu.cn}
\\Institute of Mathematics\\Nanjing Normal University\\Nanjing, 210046\\P.R.China}\end{center}\end{minipage}}}
\date{}
\maketitle
\date{}
\maketitle

\begin{abstract}
Fuzziness and randomicity widespread exist in natural science, engineering,
technology and social science. The purpose of this paper is to present a new logic - uncertain propositional logic which can deal with both fuzziness by taking truth value semantics and randomicity  by taking probabilistic semantics or possibility semantics.  As the first step for purpose of establishing a logic
system which completely reflect the uncertainty of the objective
world, this logic will lead to a set of logical foundations for
uncertainty theory as what classical logic done in certain or
definite situations or circumstances.


\end{abstract}

\parindent 24pt

\section{Introduction}

As one of the most important and one of the most widely used
concepts in the whole area of modern academic or technologic
research, uncertainty has been involved into study and applications
more than twenty years. Now along with the quickly expanding
requirements of developments of science and technology, people are
having to face and deal with more and more problems tangled with
uncertainty in the fields of natural science, engineering or
technology or even in social science. To these uncertain problems,
many traditional theories and methods based on certain conditions or
certain circumstances are not so effective and powerful as them in
the past. So the importance of research on uncertainty is emerging
more and more obviously and imminently.

We consider two well known types of uncertainty. First is fuzziness.
What is fuzziness?  from the point view of logic, a proposition $P$
is said to be a fuzzy proposition if $P$ dose not satisfy the law of
excluded middle, i.e. the formula $P \vee \neg P$ is not a
tautology. Examples of fuzzy propositions are numerous, for example,  the proposition ``This winter is cold" is a fuzzy proposition.  The second is randomicity. From the point view of logic,
randomicity can be explained as the law of causality does not hold,
i.e. the formula $(P \& (P \Rightarrow Q))\Rightarrow Q $ is not a
tautology. For example,  ``If a person fell to the ground  from 300 feet high places, he will be died. Peter fell to the ground from a 300 feet high building, he will be died". We know whether Peter will be died or not is based on various objective conditions including Peter's physical condition, the ground condition and other conditions like gravity. This reminds us to consider the following five laws in
classical logic:

(FL1) The Law of Identity,

(FL2) The Law of Contradiction,

(FL3) The Law of Excluded Middle,

(FL4) The Law of Causality,

(FL5) The Law of Sufficient Reason.

If we deny one or more of the five laws, uncertainty will occur. In
other words, uncertainty can be divided into the following five
basic types, each one corresponding to an impairment of a basic law
of formal logic:

(UCT1) \makebox[2.7cm][l]{Instability} $\sim$ Impairment of ``The Law of Identity",

(UCT2) \makebox[2.7cm][l]{Inconsistency} $\sim$ Impairment of ``The Law of Contradiction",

(UCT3) \makebox[2.7cm][l]{Fuzziness} $\sim$ Impairment of ``The Law of Excluded Middle",

(UCT4) \makebox[2.7cm][l]{Randomicity} $\sim$ Impairment of ``The Law of Causality",

(UCT5) \makebox[2.7cm][l]{Incompleteness} $\sim$ Impairment of ``The
Law of Sufficient Reason".

\vskip 3mm

We can describe these five basic types of uncertainty more elaborately as follows:

\vskip 2mm

\leftskip 5.85cm\parindent -5.05cm
(UCT1$^{\prime}$) \makebox[2.9cm][l]{Instability}\makebox[4.2mm][l]{$\sim$}Relatively\,\ ``certain"\,\ of\,\ process\,\ with ``uncertain"\,\ of properties,

(UCT2$^{\prime}$) \makebox[2.9cm][l]{Inconsistency}\makebox[4.2mm][l]{$\sim$}Relatively\,\ ``certain"\,\ of\,\ factors\,\ with\,\ ``uncertain"\,\ of\,\ relations,

(UCT3$^{\prime}$) \makebox[2.9cm][l]{Fuzziness}\makebox[4.2mm][l]{$\sim$}Relatively ``certain" of property with ``uncertain" of confines,

(UCT4$^{\prime}$) \makebox[2.9cm][l]{Randomicity}\makebox[4.2mm][l]{$\sim$}Relatively ``certain" of causation with ``uncertain" of results,

(UCT5$^{\prime}$) \makebox[2.9cm][l]{Incompleteness}\makebox[4.2mm][l]{$\sim$}Relatively\,\ ``certain"\,\ of\,\ results\,\ with\,\ ``uncertain"\,\ of\,\ causations.

\leftskip 0cm\parindent 24pt

\vskip 3mm

We know that logic is a reflection of the most basic laws of objective world in our consciousness. But, since the need of
effective thinking, this reflection is only a definite skeleton - a veriest abstract - of the various concrete representations of those
basic laws of the mutative factual world and hence there are various disparities between logic and the factual world even only between
logic and those most basic laws.

In a circumstance that people need only to consider or process certain or definite objects, this disparity will not cause any
disturbance since the disparity has been ignored in this step of consideration or processing; but if we have to face and process uncertainty
itself, then it cannot be ignored anymore, this disparity becomes serious and obvious. Therefore, when uncertainty becomes
as our main object to be considered, thought and processed, besides the logic for thinking, people must also consider the concrete basic
laws themselves of the factual world which produce the abstract basic laws of logic in our consciousness, and the relation between these two
sides.

There are various logics in the literature, weaker than classical logic, that
have been proposed to account for the failure of one of the five laws. Typical examples of them are fuzzy logic formed for reasoning fuzziness and probability logic formed for reasoning under stochastic uncertainty. There is clear distinction between fuzzy logic and probability logic, that is fuzzy logics are truth-functional while  probability logics are not  truth-functional.

Started by the pioneer paper of  L. Zadeh [19], fuzzy logic has been investigated extensively and deeply(see [6], [7], [9], [10], [23]). Now it is widely accepted that fuzzy logic is a branch of mathematical logic.  On the other hand, based on probability theory,  probability logic has been developed as an important tool for uncertainty reason (see [1], [2], [18]).

  In [18],  Nilsson presented a generalization of the ordinary true-false semantics for logical sentences to a semantics that allows probabilistic values on sentences, and described a method for probabilistic entailment in a way that probabilistic logical entailment reduces to ordinary logical
entailment when the probabilities of all sentences are either $0$ or $1$.  We will show in section 5 that every probability function $P$ can be extended to an  $\mathfrak{L}_{[0,\infty]}$-evaluation on the universe generated by those given events, this means that Nilsson's  probabilistic logic is a special case of GPL that the truth values of sentences are probability values (between $0$ and $1$).  In[15], Liu  presented an uncertain propositional logic in
which every proposition is abstracted into a Boolean uncertain variable and the truth value is defined as the uncertain measure(a generalization of probability function) that the proposition is true.  Liu's uncertain propositional logic can deal with uncertainty in the sense of Liu,  but it assumed the ``The Law of Excluded Middle", so can't deal with propositions which contain fuzziness.

Consider the following two sentences

 1.  ``There will be heavy rain at 12.00 clock tomorrow";

 2.   ``This virus will cause very serious damage on human being in the near future".

 These two sentences are neither completely random events nor fuzzy sentences since they contain both randomness and fuzziness.  So it can't be quantified independently either by probability measure or by possibility measure. This means that it is not good if we use probabilistic logic or Liu's uncertain logic to deal with these sentences. So we need a new logic system which can reasoning under both fuzziness and randomicity.

 On the other hand, from the point view of mathematics, to dealing with mixed probabilistic/non-probabilistic
uncertainty, the theory of fuzzy probability theory has been developed. It was initiated with the introduction of fuzzy random variables by Kwakernaak (see [13], [14])
in 1978/79. Subsequently, a series of extensions
and generalizations have been reported. The developments show differences in
terminology, concepts, and in the associated consideration
of measurability, and the investigations are ongoing. For an
overview on the developments and unifying considerations, see
[8] and [12]. In fuzzy probability theory, since the data from which the probabilities must be estimated are usually
incomplete, imprecise or not totally reliable, the probabilities are assumed to be fuzzy
rather than real numbers.

 The purpose of this paper is to introduce a new propositional logic which contains fuzzy logic as one of its theories and also can cope with uncertainty due to randomicity by taking probabilistic semantics or possibility semantics.

 \section{UL-algebras}

 For the purpose of dealing with semantics of uncertain propositional logic, we need to introduce a
 special algebra, which we call UL-algebra.

\begin{definition} A UL-algebra is an algebra $ \mathfrak{L} = (L, \leq, \ast, \oplus, \neg, 0, 1)$
 with two binary operations $``\ast",\  ``\oplus",\ $ and an unary operation $``\neg"$  and two
 constants $0, 1 \in L$ such that

(U1) $(L, \leq)$ is a poset with the least element $0$;

(U2) $(L, \ast, 1)$ is a commutative monoid with the unit element
$1$ and  the operation $\ast$ can be restricted to the interval $[0, 1] = \{a\in L \mid a \leq 1\}$ such that it becomes a sub-monoid;

(U3)  $(L, \oplus, 0)$ is a commutative monoid with the unit
element $0$.  If $x \leq y$ in $L$ then $x
\oplus z \leq y \oplus z$ for all $z\in L$. For $a, b \in L$ and $a \leq b$, there exists $max\{c \in L\mid a \oplus c = b\}$;

(U4) $\neg: L \rightarrow L$ is an unary operation such that $\neg 0 =1, \neg 1 = 0$,  and $\neg b\leq \neg a$ whenever $a \leq b$;

(U5)For every  $a \in L$,  the mapping  $a\ast ( ): L \rightarrow L, x \mapsto a \ast x$ has a right adjoint $a \rightarrow ( )$,  i.e.  $y \leq a
\rightarrow x$ if and only if $a \ast y \leq x$ for all $x, y \in L$.

 \end{definition}

\begin{example} 1. Let $L$ be a Heyting algebra. Then $L$ is an UL-algebra when we define $a \ast b = a \wedge b$ be the meet for $\{a, b\}$,  $a \oplus b= a \vee b$ be the join for $\{a, b\}$,  $\neg a= a \rightarrow 0$. More generally,  every BL-algebra or MTL-algebra is an UL-algebra.

2. Consider the right interval $[0, \infty] = \{x \in \mathbb{R}\mid x \geq 0\} \cup \{\infty\}$. Define an order on $[0, \infty]$ as following:

For $x , y \in \{x \in \mathbb{R}\mid x \geq 0\}$,  $x \leq y$ in the usual sense and $x \leq \infty$ for all $x \in  [0, \infty]$.

Consider the following two types of $\ast$ operation

(i)  $x \ast y= min(x, y)$;

(ii)  $x \ast^{\prime} y = x \times y$ for $x, y\not= \infty$,   $0 \ast^{\prime} \infty = \infty \ast^{\prime} 0 = 0$,  $a \ast^{\prime} \infty = \infty \ast^{\prime} a = \infty$ for $a \not= 0$.

Define

(i $^{\prime}$)  $x \oplus y = x \vee y$;

 (ii $^{\prime}$) $x \oplus^{\prime} y = x + y$ for $x, y \not= \infty$ and $x \oplus^{\prime} \infty = \infty \oplus^{\prime} x = \infty$ for all $x \in [0, \infty]$.

  Define $\neg x = 1- x$ for $x \leq 1$ and $\neg x = 0$ for all $x \not\leq 1$.

  Then $([0, \infty], \leq, \ast, \oplus, \neg, 0, \infty)$, $([0, \infty], \leq, \ast, \oplus^{\prime}, \neg, 0, \infty)$,  $([0, \infty], \leq, \ast^{\prime}, \oplus, \neg, 0, 1)$  and  $([0, \infty], \leq, \ast^{\prime}, \oplus^{\prime}, \neg, 0, 1)$ are all UL-algebras.

  3. Let $\ast$ be a continuous t-norm on $[0, 1]$. Define $x \oplus y = x \vee y, \neg x = 1 - x$, then $([0, 1], \leq, \ast, \oplus, \neg, 0, 1)$ is a UL-algebra.
\end{example}

\begin{lemma} Let $\mathfrak{L} = (L, \leq, \ast, \oplus, \neg, 0, 1)$ be an UL-algebra, and $x, y, z \in L$. The
following properties hold:

(1) If $x \leq y$, then $x \ast z \leq y \ast z, \, (z\rightarrow x)
\leq (z\rightarrow y),\, (y\rightarrow z) \leq (x\rightarrow z)$;

(2) $ x \leq x\oplus y$.

(3) If $x \rightarrow y = 1$ then $x \leq y$. The converse is not true.

(4) $(1 \rightarrow x) = x$.
\end{lemma}

{\bf Proof}: (1) Suppose $x \leq y$. We have $y \leq z\rightarrow
(y \ast z)$ since $(y \ast z) \leq (y \ast z)$. Hence $x \leq
(z\rightarrow (y \ast z))$, this implies $(x \ast z) \leq (y \ast
z)$.  $(z\rightarrow
x) \leq (z\rightarrow y)$ since $z \ast (z\rightarrow x) \leq x \leq y$;\,  $(y\rightarrow z) \leq (x\rightarrow z)$ since $x \ast (y\rightarrow z)
 \leq  y \ast (y\rightarrow z) \leq z$.

(2)  $x = x
\oplus 0 \leq x \oplus y$.

(3) $x =x \ast 1 = x \ast (x \rightarrow y ) \leq y$. Take an UL-algebra $\mathfrak{L} = ([0, \infty], \leq, \ast^{\prime}, \oplus, \neg, 0, 1)$ be the UL-algebra  given in Example 2.1. Take $x = 0$ and $y = 0.2$ then we have $x \rightarrow y = \infty$.

(4) Clear.
$\hfill\Box$

\begin{lemma} Let $\{\mathfrak{L} = (L_{i}, \leq_{i}, \ast_{i}, \oplus_{i}, \neg_{i}, 0_{i}, 1_{i}) \mid i\in I\}$ be a family of UL-algebras. The product
$\prod_{i\in I}L_{i}$ of $L_{i}$ with pointwise order and pointwise operations is a UL-algebra which is denoted by $\prod_{i\in I}\mathfrak{L}$. \end{lemma}

\section{The Basic Uncertain Propositional Logic}

In the basic uncertain propositional logic, similar to fuzzy logic and probabilistic logic, we will start with some atomic
formulas which are basic materials to construct any complicated
uncertain proposition. We will not distinguish the words
``formula",
 ``proposition" and  ``sentence".

The language of uncertain
propositional logic contains atomic formulas denoted by $\varphi,
\psi, \cdot\cdot\cdot,$ and a constant  $\bar{0}$ which we call
false; words or phrases for ``and",  ``or",  ``not",  ``implies"
and  ``if and only if" are presented by
 symbols as follows:

 $\&$ for  ``and", $\vee$ for  ``or", $\neg$ for  ``not", $\Rightarrow$ for  ``implies", $\Leftrightarrow$
 for ``if and only if".

 Let $\varphi$ and $\psi$ be two sentences. Then

 $\neg\varphi$ means ``not $\varphi$";

 $\varphi \& \psi$ means ``$\varphi$ and $\psi$";

 $\varphi \vee \psi$ means ``$\varphi$ or $\psi$";

  $\varphi \Rightarrow \psi$ means ``if $\varphi$ then $\psi$";

   $\varphi \Leftrightarrow \psi$ means ``$\varphi$ if and only if $\psi$".

Propositions are defined in the obvious way: each atomic formula
 is a proposition; $\bar{0}$ is a proposition; if $\varphi$ and $\psi$ are
 propositions, then $\varphi \& \psi, \varphi  \vee \psi, \varphi \Rightarrow
 \psi, \neg \varphi$ are propositions. Further
 connectives are defined as follows:

 $\varphi \Leftrightarrow \psi$ is an abbreviation for  $ (\varphi \Rightarrow \psi)\& (\psi
 \Rightarrow \varphi)$;

$\phi\mid\psi$ is an abbreviation for $\psi \Rightarrow (\psi \&\phi)$.

We call formula $\phi\mid\psi$ the condition formula of $\phi$ under $\psi$.

Consider the following formulas:

(A1) $(\varphi \Rightarrow \psi) \Rightarrow ((\psi \Rightarrow \chi) \Rightarrow (\varphi \Rightarrow
 \chi))$;

 (A2)  $\varphi \Rightarrow \varphi$;

(A3) $(\varphi \Rightarrow (\psi\Rightarrow \chi)) \Rightarrow (\psi \Rightarrow (\varphi\Rightarrow \chi))$;

 (A4) $(\varphi \& \psi) \Rightarrow \varphi$;

 (A5) $(\varphi \& \psi) \Rightarrow (\psi \& \varphi)$;

(A6)   $(\varphi\vee \psi) \Rightarrow (\psi\vee
  \varphi),  (\bar{0} \vee
  \varphi)\Rightarrow \varphi$;

(A7) $\bar{0} \Rightarrow \varphi$.

\begin{definition} The formulas A1-A7 are axioms of the basic uncertain propositional logic(shortly UPL).
  A proof in UPL is a sequence $\varphi_{1}, \cdot\cdot\cdot, \varphi_{n}$ of formulas such
  that each $\varphi_{i}$  is either an axiom of UPL or follows from
  some preceding $\varphi_{j}, \varphi_{k} (j,k < i)$ by modus ponens: from $\varphi$ and $\varphi
  \Rightarrow \psi$ infer $\psi$. A formula is provable in UPL(notation: $\vdash \varphi$)
  if it is the last member of a proof.\end{definition}

 Let $\mathfrak{L} = (L, \leq, \ast, \oplus, \neg,  0, 1)$ be an UL-algebra. An $\mathfrak{L}$-evaluation (or an $\mathfrak{L}$-model) is a mapping $e$ from the set of all atomic formulas into $[0, 1] = \{a \in L \mid a \leq 1\}$ with an extension  from the set of all formulas into $L$(which we also denoted as $e$) satisfying the following conditions:

 $e(\bar{0}) = 0$;

 $e( \varphi \Rightarrow \psi) =  e(\varphi) \rightarrow e(\psi)$;

 $e(\neg\varphi)  = \neg e(\varphi)$;

 $e( \varphi \& \psi) \leq e(\varphi),   e( \varphi \& \psi) \leq  e(\psi), e( \varphi \& \psi) = e(\psi \& \varphi)$,  $e( \varphi \& \psi) = e(\psi)$ if $1 \leq e(\varphi)$;

 $e( \varphi \vee \psi) =  max\{c \in L\mid e(\varphi \& \psi) \oplus c = e(\varphi)\oplus e(\psi)\}$ ( written  $ e(\varphi) \overline{\oplus} e(\psi))$.

Note that for every $\mathfrak{L}$-evaluation $e$,   $e(\varphi \vee \psi) =   e(\psi \vee \varphi)$ since $e(\varphi \& \psi) =   e(\psi \& \varphi)$.

 \begin{definition} For a given UL-algebra $\mathfrak{L} = (L, \leq, \ast, \oplus, \neg,  0, 1)$, a formula $\varphi$ is called a $\mathfrak{L}$-tautology
 if for every $\mathfrak{L}$-evaluation $e$,  $1 \leq e(\varphi)$.  $\varphi$ is called a tautology if $\varphi$ is a $\mathfrak{L}$-tautology for every UL-algebra $\mathfrak{L}$.
 \end{definition}

 \begin{proposition} The formulas A1-A7 are tautologies.
\end{proposition}

  {\bf Proof}: (A1) Let $\mathfrak{L} = (L, \leq,  \ast,  \oplus,  \neg,  0,  1)$ be an UL-algebra. It suffice to show $1 \leq (a \rightarrow b) \rightarrow ((b\rightarrow c) \rightarrow (a \rightarrow c)) $ for all $a,  b,  c \in L$, but it is clear.

  (A2) Clear.

  (A3) Let $\mathfrak{L} = (L, \leq, \ast,  \oplus,  \neg,  0, 1)$ be an UL-algebra. It suffice to show $1 \leq (a \rightarrow (b\rightarrow c)) \rightarrow (b \rightarrow (a\rightarrow c))$ for all $a, b, c \in L$, this is clear.

  (A4) Let $\mathfrak{L} = (L, \leq, \ast, \oplus, \neg,  0, 1)$ be an UL-algebra and $e$ an $\mathfrak{L}$-evaluation. We have $e(\varphi \& \psi) \leq e(\varphi)$.  Hence $1 \leq e((\varphi \& \psi) \rightarrow \varphi)$.

  (A5) Clear.

  (A6)  $(\varphi\vee \psi) \Rightarrow (\psi\vee
  \varphi)$ is clear. For an UL-algebra $\mathfrak{L} = (L, \leq, \ast, \oplus, \neg,  0, 1)$ and an $\mathfrak{L}$-evaluation $e$, we have $e(\bar{0} \vee
  \varphi) = e(\varphi)$ since $e(\bar{0}\&
  \varphi) \leq e(\bar{0}) = 0$.

  (A7) Clear.  $\hfill\Box$

  We write $\bar{1}$ for $\bar{0} \Rightarrow \bar{0}$.

  \begin{proposition} The following properties are  provable in UPL:

  (1) $\bar{1}$;

  (2) $(\bar{1} \Rightarrow \varphi) \Rightarrow \varphi$;

  (3) $\varphi \Rightarrow (\varphi \&(\psi\mid \varphi) \Rightarrow \psi)$.
 \end{proposition}

{\bf Proof}:  (1) Clear by the definition of $\bar{1}$ and (A2).

 (2)  $UPL \vdash (\bar{1} \Rightarrow \varphi) \Rightarrow (\bar{1} \Rightarrow \varphi)$ by (A2), and

$ UPL \vdash ((\bar{1} \Rightarrow \varphi) \Rightarrow (\bar{1} \Rightarrow \varphi)) \Rightarrow (\bar{1} \Rightarrow ((\bar{1} \Rightarrow \varphi) \Rightarrow \varphi))$ by (A3). Thus

$ UPL \vdash \bar{1} \Rightarrow ((\bar{1} \Rightarrow \varphi) \Rightarrow \varphi)$ by modus ponens.

  $UPL \vdash (\bar{1} \Rightarrow \varphi) \Rightarrow \varphi$ by (1) and modus ponens.

  (3) $UPL \vdash \varphi \& (\psi\mid\varphi) \Rightarrow (\varphi \Rightarrow \varphi \& \psi)$ by the definition of $(\psi\mid\varphi)$ and (A4). Thus

  $ UPL \vdash \varphi \Rightarrow (\varphi \&(\psi\mid \varphi) \Rightarrow \psi)$ by (A3)  $\hfill\Box$

  \begin{theorem} The logic UPL is sound with respect to tautologies: if $\varphi$ is provable in UPL, then
  $\varphi$ is a tautology.\end{theorem}

  {\bf Proof} We have shown that all axioms of UPL are
  tautologies.  For the deduction rule of UPL:
   suppose $\phi$ and $\phi \Rightarrow \psi$ are both tautologies, then for every UL-algebra $\mathfrak{L} = (L, \leq, \ast, \oplus, \neg, \rightarrow, 0, 1)$ and an $\mathfrak{L}$-evaluation $e$, we have $1 \leq e(\varphi) \ast e(\varphi \rightarrow \psi) = e(\varphi) \ast (e(\varphi) \rightarrow e(\psi)) \leq e(\psi)$. Hence $\psi$  is a tautology. $\hfill\Box$

 \begin{definition} For a given UL-algebra $ \mathfrak{L} = (L, \leq, \ast, \oplus, \neg,  0, 1)$, consider some $\mathfrak{L}$-tautologies $\mathcal{L}$. We will call the collection of propositions A(1)- A(7) plus $\mathcal{L}$  an $\mathfrak{L}$  uncertain theory over UPL (or $\mathfrak{L}$ logic). A
 proof in an $\mathfrak{L}$ uncertain theory is a sequence of formulas: $\varphi_{1},
 \cdot\cdot\cdot, \varphi_{n}$ where each $\varphi_{i}$ is either a member of the  uncertain theory or follows from some preceding
 member of the sequence by modus ponens. We write $\mathfrak{L} \vdash \varphi$
 for $\varphi$ to be provable in the $\mathfrak{L}$ uncertain theory, i.e. the last member
 of a proof in the $\mathfrak{L}$ uncertain theory.\end{definition}

 We call an  $\mathfrak{L}$  uncertain theory  inconsistent if $\mathfrak{L} \vdash \bar{0}$, otherwise
 it is consistent. By A(8), we have the following result.

 \begin{lemma} An  $\mathfrak{L}$ uncertain theory is inconsistent if and only if $\mathfrak{L} \vdash
 \psi$ for each $\psi$.\end{lemma}

  \begin{example} Consider the case for UL-algebra being the the two-elements lattice $2 = \{0, 1\}$ with $\ast = \wedge,\,  \oplus = \vee,\, \neg = 1- ( )$. Then we can add some axioms(all classical laws) such that the $2$ logic is just the classical propositional logic, i.e. the classical Boolean logic is a theory over UPL. \end{example}

 \section{Generalized Fuzzy Logic}

 In this section we consider the uncertain theory over UPL for  UL-algebra being $\mathfrak{I}_{f} = (I,  \leq,  \ast,   \oplus,   \neg,  0,  1)$ where $I = [0, 1]$ be the unit interval, $\leq$ be the natural order and

 $x \ast y = x\wedge y$ be the join of $x$ and $y$;

  $x \oplus y = x \vee y$ be the meet of $x$ and $y$;

    $\neg x = 1- x$ for all $x, y \in I$.

    Under this case, an $\mathfrak{I}_{f}$-evaluation $e$ will satisfy

    $e(\varphi  \& \psi) \leq e(\varphi) \wedge e(\psi)$;

    $ e(\varphi \vee \psi) = e(\varphi) \vee e(\psi)$.

    We call this theory generalized fuzzy logic, written as GFL.

 \begin{definition} The axioms of GFL are those of UPL plus

 (GFL1) $((\varphi \Rightarrow \psi) \Rightarrow \chi)) \Rightarrow (((\psi \Rightarrow \varphi) \Rightarrow \chi) \Rightarrow \chi)$;

 (GFL2)  $(\varphi \Rightarrow (\psi \Rightarrow \chi)) \Rightarrow ((\varphi \& \psi)
 \Rightarrow \chi)$;

 (GFL3)  $\varphi \Rightarrow (\psi \Rightarrow \varphi)$;

 (GFL4)  $\varphi\Rightarrow (\varphi \vee \psi)$;

 (GFL5) $(\varphi \vee \psi) \vee \chi \Leftrightarrow \varphi \vee (\psi \vee \chi)$;

 (GFL6) $(\varphi \vee \psi) \Rightarrow (\varphi \Rightarrow \psi) \Rightarrow \psi)$.

 (GFL7) $(\neg\neg \varphi) \Leftrightarrow \varphi$.
\end{definition}

\begin{lemma}  (GFL1)-(GFL6) are all
 $\mathfrak{I}_{f}$-tautologies.
\end{lemma}

{\bf Proof}  (GFL1) Let $e$ be an  $\mathfrak{I}_{f}$-evaluation. Then we have $e(\varphi) \leq e(\psi)$ or $ e(\psi) \leq e(\varphi)$ holds. Under both cases we have $e(((\varphi \Rightarrow \psi) \Rightarrow \chi)) \Rightarrow (((\psi \Rightarrow \varphi) \Rightarrow \chi) \Rightarrow \chi)) =1$.

(GFL2) We only need to note that for any  $\mathfrak{I}_{f}$-evaluation $e$, we have $e(\varphi \& \psi) \leq e(\varphi) \wedge e(\psi)$.

(GFL3) For any  $\mathfrak{I}_{f}$-evaluation $e$, we have $e(\varphi) \wedge e(\psi) \leq e(\varphi)$. Hence $e(\varphi \Rightarrow (\psi \Rightarrow \varphi)) = 1$.

(GFL4)-(GFL5) are clear by the fact that for any $\mathfrak{I}_{f}$-evaluation $e$ we have $e(\varphi \vee \psi) = e(\varphi) \vee e(\psi)$.

(GFL6) For a given $\mathfrak{I}_{f}$-evaluation $e$, if $e(\varphi) \leq e(\psi)$ then $e((\varphi \vee \psi) \Rightarrow (\varphi \Rightarrow \psi) \Rightarrow \psi)) = e(\psi) \rightarrow e(\psi) = 1$; If $e(\varphi) \not\leq e(\psi)$ then $e((\varphi \vee \psi) \Rightarrow (\varphi \Rightarrow \psi) \Rightarrow \psi)) = e(\varphi) \rightarrow (e(\psi)\rightarrow e(\psi)) = 1$.

(GFL7) Clear.
$\hfill\Box$

\begin{corollary} Every provable proposition $\varphi$ in GFL is an $\mathfrak{I}_{f}$-tautology.
\end{corollary}

\begin{proposition}  The following formulas are provable in GFL:

(1) $(\varphi \& (\varphi \Rightarrow \psi)) \Rightarrow \psi$, in particular, $\varphi \& \psi\mid\varphi \Rightarrow \varphi \& \psi$;

(2) $(\varphi \Rightarrow \psi) \vee (\psi \Rightarrow \varphi)$;

(3) $\neg\neg\varphi \Rightarrow (\psi \Rightarrow \varphi)$;

(4) $(\varphi \Rightarrow \psi) \Rightarrow ((\varphi \vee \psi) \Rightarrow \psi)$.
\end{proposition}

{\bf Proof} (1) $ GFL\vdash (\varphi \Rightarrow \psi)
\Rightarrow (\varphi \Rightarrow \psi)$, hence

$ GFL\vdash (\varphi \Rightarrow ((\varphi\Rightarrow\psi) \Rightarrow \psi)$ by (A3).

$ GFL\vdash (\varphi \& (\varphi \Rightarrow \psi)) \Rightarrow \psi$ by (GFL2).

(2) $GFL\vdash (\varphi \Rightarrow \psi) \Rightarrow (\varphi \Rightarrow \psi) \vee (\psi \Rightarrow \varphi)$ and

$ GFL\vdash (\psi \Rightarrow \varphi) \Rightarrow (\varphi \Rightarrow \psi) \vee (\psi \Rightarrow \varphi)$ by (GFL4). Thus

$ GFL\vdash (\varphi \Rightarrow \psi) \vee (\psi \Rightarrow \varphi)$ by (GFL1) and modus ponens.

(3) $GFL\vdash \neg\neg\varphi \Rightarrow \varphi$ by (GFL7) and (A4). Hence

 $GFL\vdash \neg\neg\varphi \Rightarrow (\psi \Rightarrow \varphi)$ by (GFL3) and (A1).

 (4) $GFL\vdash (\varphi \vee \psi) \Rightarrow (\varphi \Rightarrow \psi) \Rightarrow \psi)$;

 Hence  $GFL\vdash (\varphi \Rightarrow \psi) \Rightarrow((\varphi \vee \psi) \Rightarrow \psi)$ by (A3).
 $\hfill\Box$

 Note though  that the logic GFL seems similar to
Godel fuzzy logic,  but each of them is not a theory of another. For example,  (GFL7) is not a theorem of Godel logic and the axiom $\varphi \Rightarrow (\varphi \& \varphi)$ of Godel logic is not a tautology in GFL.

\section{Generalized Probability Logic}

 In this section we shall consider the uncertain theory over UPL for UL-algebra being $ \mathfrak{L}_{[0,\infty]} = ([0, \infty],  \leq,  \times,  +,  \neg,   0,  1)$ where $[0, \infty] = \{x \in \mathbb{R}\mid x \geq 0\} \cup \{\infty\}$,  for $x,  y \in \{x \in \mathbb{R}\mid x \geq 0\}$,

 $x \leq y$ in the usual sense  and $x \leq \infty$ for all $x \in  [0, \infty]$;

 $x \ast y = x \times y$ for $x, y \not= \infty$,   $0 \ast \infty = \infty \ast 0 = 0$,  $a \ast \infty = \infty \ast a = \infty$ for $a \not= 0$;

  $x \oplus y = x + y$ for $x, y \not= \infty$ and  $x \oplus \infty = \infty \oplus x = \infty$ for all $x \in [0, \infty]$;

   $\neg x = 1- x$ for $x \leq 1$ and $\neg x = 0$ for all $x \not\leq 1$.

  Under this case, an $\mathfrak{L}_{[0,\infty]}$-evaluation $e$ will satisfies

  $e(\varphi  \& \psi) \leq e(\varphi) \wedge e(\psi)$ and

   if $e(\varphi), e(\psi) \not= \infty$ then $ e(\varphi  \vee \psi) = e(\varphi) + e(\psi) - e(\varphi \& \psi)$,

   if $e(\varphi)  = \infty$, or $ e(\psi) = \infty$ then $ e(\varphi  \vee \psi) = \infty$.

   We call this theory generalized probability logic, written as GPL.

 The following formulas GP1-GP3 plus A1-A7 are axioms of the generalized probability  logic GPL.

 (GPL1) $ \neg\neg\varphi \Rightarrow \varphi$;

 (GPL2)  $\varphi\Rightarrow (\varphi \vee \psi)$;

 (GPL3) $(\varphi \& \neg\varphi) \Leftrightarrow \neg(\varphi \vee \neg\varphi)$.

 \begin{lemma}  (GPL1)-(GPL3) are all $\mathfrak{L}_{[0,\infty]}$-tautologies. \end{lemma}

{\bf Proof} (GPL1) is clear.

 (GPL2) Let $e$ be an  $\mathfrak{L}_{[0,\infty]}$-evaluation. If $e(\varphi) = \infty$ or $e(\psi) = \infty$ , then $e(\varphi\Rightarrow (\varphi \vee \psi)) = \infty$; if $e(\varphi), e(\psi) \not= \infty$, then $e(\varphi  \vee \psi) = e(\varphi) + e(\psi) - e(\varphi \& \psi) \geq e(\varphi)$ since $e(\varphi \& \psi) \leq e(\psi)$.

  (GPL3) For every  $\mathfrak{L}_{[0,\infty]}$-evaluation $e$,  if $e(\varphi) \not\leq 1$ then $e(\neg\varphi) = 0,   e(\varphi \& \neg\varphi) = 0,  e(\varphi \vee \neg\varphi) = e(\varphi)$,  hence
  $e(\varphi \& \neg\varphi) = e(\neg(\varphi \vee \neg\varphi)) = 0$;  if $e(\varphi) \leq 1$ then $e(\varphi \vee \neg\varphi) =  e(\varphi) + e(\neg\varphi) - e(\varphi \& \neg\varphi) = 1 -  e(\varphi \& \neg\varphi)$, i.e.  $e(\varphi \& \neg\varphi) = e(\neg(\varphi \& \neg\varphi))$.  This shows $1 \leq e((\varphi \& \neg\varphi) \Leftrightarrow \neg(\varphi \vee \neg\varphi))$.
$\hfill\Box$

\begin{corollary} Every provable proposition $\varphi$ in GPL is an $\mathfrak{L}_{[0,\infty]}$-tautology.
\end{corollary}

\begin{proposition}   The following formulas are  provable in GPL:

(1) $\varphi \& (\psi\mid\varphi) \Rightarrow (\neg\neg\varphi \Rightarrow \varphi)$;

(2)  $\neg(\varphi \vee \neg\varphi) \Rightarrow \varphi$.
\end{proposition}

{\bf Proof} (1) $GPL\vdash \neg\neg\varphi \Rightarrow (\varphi \&(\psi\mid \varphi) \Rightarrow \psi)$ by Proposition 3.2(3) and (GPL1). Thus

$GPL\vdash \varphi \& (\psi\mid\varphi) \Rightarrow (\neg\neg\varphi \Rightarrow \varphi)$ by (A3) and modus ponens.

(2) $GPL\vdash \neg(\varphi \vee \neg\varphi) \Rightarrow (\varphi \& \neg\varphi)$ by (GPL3),

$GPL\vdash (\varphi \& \neg\varphi) \Rightarrow \varphi$ by (A4), hence

$GPL\vdash \neg(\varphi \vee \neg\varphi) \Rightarrow \varphi$ by (A1).
 $\hfill\Box$

 Recall the Kolmogorov's classical definition of probability:

 Let $\Omega$ be a set(the ``universal set"), and let $\mathfrak{F}$ be a field on $\Omega$, i.e. $\mathfrak{F}$ be a set of subsets of $\Omega$ which contains
 $\Omega$ as a member and closed under complementation and finite union. Members in $\mathfrak{F}$ are called events. A probability function on $\mathfrak{F}$ is
 a function $P: \mathfrak{F} \rightarrow [0, 1]$ satisfying conditions:

 (P1) $P(\Omega) = 1$;

 (P2)  $P(A \cup B) = P(A) + P(B)$ for all $A, B \in \mathfrak{F}$ such that $A \cap B = \emptyset$.

 Call such a triple of $(\Omega, \mathfrak{F}, P)$ a probability space.

 By classical results we know that the conditions (P1) and (P2) in the above definition can be equivalently replaced by the following three conditions:

 (P1$^{\prime}$)  $P(\emptyset) = 0$;

 (P2$^{\prime}$)  $P(A^{c}) = 1 - P(A)$;

 (P3$^{\prime}$)  $P(A \cup B) = P(A) + P(B) - P(A\cap B)$.

 We expand the collection of events $\mathfrak{F}$ to a new collection of propositions  $\widetilde{\mathfrak{F}}$ by add a new connective ``$\Rightarrow$" as following:

 Each member of $\mathfrak{F}$ is a member of $\widetilde{\mathfrak{F}}$;  if $\varphi$ and $\psi$ are
members of $\widetilde{\mathfrak{F}}$,  then $\varphi \& \psi, \varphi  \vee \psi, \varphi \Rightarrow
 \psi, \neg \varphi$ are members of $\widetilde{\mathfrak{F}}$. Here if  $\varphi$ and $\psi$ are
members of $\mathfrak{F}$, we define $\varphi \& \psi$ is same as $ \varphi \cap \psi$,  $\varphi  \vee \psi$ is same as $\varphi \cup \psi$ and $\neg \varphi $ is same as
$\varphi^{c}$.

 \begin{theorem} Let $(\Omega, \mathfrak{F}, P)$ be a probability space. Then the probability function $P$ can be extended to an  $\mathfrak{L}_{[0,\infty]}$-evaluation on  $\widetilde{\mathfrak{F}}$.
 Conversely,  given an  $\mathfrak{L}_{[0,\infty]}$-evaluation $e$ on  $\widetilde{\mathfrak{F}}$ such that $e(\varphi) \leq 1$ for each member $\varphi$ of $\mathfrak{F}$. The restriction of $e$ to  $\mathfrak{F}$ is a probability function on
 $\mathfrak{F}$. \end{theorem}

 \section{Fuzzy Random Logic}

 We know that fuzzy logics only consider fuzziness and do not contains randomicity or other uncertainty while probability logic just reasoning under randomicity  and contains no fuzziness. But in the real world, fuzziness and  randomicity and other uncertainties are in fact co-existed.  To deal with propositions which contain both fuzziness and randomicity we consider  the uncertain theory over UPL for  UL-algebra being $\mathfrak{I}_{f} \times \mathfrak{L}_{[0,\infty]}$ for which the poset being the product $ [0, 1] \times [0, \infty]$ with pointwise order and pointwise operations.

 For an $\mathfrak{I}_{f} \times \mathfrak{L}_{[0,\infty]}$-evaluation $e$ and any proposition $\varphi$, $e(\varphi) = (a, b)$, we may think that the first coordinate $a$ is the degree of truth of the proposition $\varphi$ and the second coordinate $b$  the degree of belief of  $\varphi$.

 For fuzzy random logic FRL, we consider the following axioms:

 (FRL1)  $\varphi\Rightarrow (\varphi \vee \psi)$;

 (FRL2) $ \neg\neg\varphi \Rightarrow \varphi$;

 By Lemma 4.1 and 5.1 we have the following result.

 \begin{lemma} (FRL1) and (FRL2) are both $\mathfrak{I}_{f} \times \mathfrak{L}_{[0,\infty]}$-tautologies. \end{lemma}

 \begin{corollary} Every provable proposition $\varphi$ in FRL is an $\mathfrak{I}_{f} \times \mathfrak{L}_{[0,\infty]}$-tautology.
\end{corollary}

\begin{example} Consider three sentences $\varphi$, $\psi$ and $\varphi\rightarrow \psi$. If we are given the probability of $\varphi$, denoted by $p(\varphi)$; the probability of $\varphi\rightarrow \psi$, denoted by $p(\varphi\rightarrow \psi)$. These values
are estimated, or they are provided by some experts, so they are not certain. Suppose we are also given the degree of truth of the probability of $\varphi$, denoted by $t(\varphi)$ and the degree of truth of the probability of $\varphi\rightarrow \psi$, denoted by $t(\varphi\rightarrow \psi)$ respectively. We consider these data be an $\mathfrak{I}_{f} \times \mathfrak{L}_{[0,\infty]}$-evaluation $e$ such that

$e(\varphi) = (p(\varphi), t(\varphi)),   e(\varphi \rightarrow \psi) = (p(\varphi \rightarrow \psi), t(\varphi \rightarrow \psi))$.

We now can give a lower bound and upper bound of the probability $p(\psi)$ of the sentence $\psi$ and its truth degree $t(\psi)$ as following

$p(\varphi)\times p(\varphi \rightarrow \psi) \leq  p(\psi)  \leq p(\varphi \rightarrow \psi)$;

$t(\varphi)\wedge t(\varphi \rightarrow \psi) \leq  t(\psi)  \leq t(\varphi \rightarrow \psi)$.   \end{example}

\section{Properties for Uncertain Propositional Logic}

In this section we consider some basic laws in
classical logic, we will see that these laws are not tautologies in UPL.

The first is the law of excluded middle: $\varphi \vee
\neg\varphi$. We know it is is a tautology in classical logic. Now
it is clear that it is not a tautology in UPL, even in GFL and GPL. For example,  taking an $\mathfrak{I}_{f}$-evaluation $e$ such that $e(\varphi) = 0.5$ then $e(\varphi \vee
\neg\varphi) = 0.5$.

The case that ``The Law of Excluded Middle is not satisfied" corresponds to fuzzy objects or phenomena, such as concepts ``good", ``bad" (even ``not good"),
``cold", ``warm" (even ``not cold"), ``young", ``old" (even ``not young"), etc. In these dual concepts or phenomena, people cannot or cannot
reasonably cut them into two or more crisp divided parts, although they have being used more than thousands of years and will
continue to be used or must continue to be used in the future. The appearance of fuzzy set theory reflects this kind of need in processing
fuzzy objects, and this kind of situations is naturally contained in our frame of uncertainty.

\vskip 4mm

The second is that of the law of contradiction $(\varphi \&
\neg\varphi) \Rightarrow \bar{0}$. We know this is a tautology in
classical logic and fuzzy logic. But it is not a
 tautology in UPL. For example,  we take an $\mathfrak{I}_{f}$-evaluation $e$ such that $e(\varphi) = 0.5$ and $e(\varphi \& \neg\varphi) = e(\varphi) \wedge e(\neg\varphi) = 0.5$. Of course we can give an uncertain theory over UPL such that $(\varphi \&
\neg\varphi) \Rightarrow \bar{0}$ is a  tautology in this theory. For example, we consider UL-algebra $(I, \leq, \wedge, \vee, ( )\rightarrow 0,  0, 1)$.

In many factual situations or circumstances, information can be obtained from objects contains contradictive parts
or factors, they often unavoidably contradict to some others in different extents. In these cases, contradictions with various
extents become into objects we have to deal with directly in our researches or operations, but cannot always be avoided or eliminated or cleared
up. Our disposal in uncertain logic for contradictions with some extents is just to fulfill this need. In fact, these contradictions caused by
various squeeze and mixture of different levels from a high dimension space into a low dimension space; this problem will
be solved in our next paper.

\vskip 4mm

 The third  is the law of causality ($\varphi \& (\varphi\Rightarrow
\psi)\Rightarrow \psi $). We know this is a tautology in classical
logic and fuzzy logic, but not a  tautology in UPL and GPL.  For example,  consider an  $\mathfrak{L}_{[0,\infty]}$-evaluation $e$ such that $e(\varphi) = 0.8, e(\psi) = 0.4$ and $e(\varphi \& (\varphi \Rightarrow \psi)) = 0.5$ then $e((\varphi \& (\varphi \Rightarrow \psi)) \Rightarrow \psi) = 0.8$.  But $\varphi \& (\varphi\Rightarrow
\psi)\Rightarrow \psi $ is a tautology in GFL by Proposition 4.1 (1).

As well known, randomness appears from impairments of the law of causality. Although ``Does there exist objectively or absolutely
randomness?" will always be an argument without an absolute answer, but the fact that we have to always face to and dispose various processes
which cannot satisfy the law of causality is still a widely accepted choice. Obviously, randomness reduces believablity or
reliability, considering their extents in a way uniform with fuzziness is just what we have done in
the frame of uncertain propositional logic.

\end{document}